\documentclass[12pt]{article}                                               

\input{amssym.def}                                                           
\input{amssym.tex}                                                           
                                                                             
\begin{document}                                                             
\title{On traces  of  Frobenius endomorphisms }

\author{Igor  ~Nikolaev
\footnote{Partially supported 
by NSERC.}}


\date{}
 \maketitle


\newtheorem{thm}{Theorem}
\newtheorem{lem}{Lemma}
\newtheorem{dfn}{Definition}
\newtheorem{rmk}{Remark}
\newtheorem{cor}{Corollary}
\newtheorem{cnj}{Conjecture}
\newtheorem{exm}{Example}


\begin{abstract}
We compute  the number of  points of  projective variety $V$ over
a  finite field in terms of invariants  of the  so-called  Serre  $C^*$-algebra of
 $V$.

\vspace{7mm}


\vspace{5mm}

\end{abstract}

\section{Introduction}
The number of solutions of a system of polynomial equations over
a finite field is an important invariant of the system and 
an old problem dating back to Gauss.    Recall that if ${\Bbb F}_q$
is a field with $q=p^r$ elements and $V({\Bbb F}_q)$ a smooth $n$-dimensional 
projective variety over ${\Bbb F}_q$,  then  one can define a zeta
function $Z(V; t):= \exp~\left(\sum_{r=1}^{\infty} |V({\Bbb F}_{q^r})|{t^r\over r}\right)$;
the function is rational,  i.e.
\begin{equation}\label{eq1}
Z(V; t)={P_1(t)P_3(t)\dots P_{2n-1}(t)\over P_0(t)P_2(t)\dots P_{2n}(t)},
\end{equation}
where   $P_0(t)=1-t$, $P_{2n}(t)=1-q^nt$ and for each $1\le i\le 2n-1$ the polynomial 
$P_i(t)\in {\Bbb Z}[t]$ can be written as  $P_i(t)=\prod_{j=1}^{deg~P_i(t)} (1-\alpha_{ij}t)$
so that   $\alpha_{ij}$ are algebraic integers with
$|\alpha_{ij}|=q^{i\over 2}$, see e.g.  [Hartshorne 1977]  \cite{H},  pp. 454-457. 
The  $P_i(t)$ can be viewed as characteristic polynomial of the Frobenius 
endomorphism $Fr_q^i$  of  the $i$-th $\ell$-adic cohomology
group  $H^i(V)$;  such an endomorphism is induced 
by the map acting on  points  of variety $V({\Bbb F}_q)$
according to the formula  $(a_1,\dots, a_n)\mapsto (a_1^q,\dots, a_n^q)$.   
(We assume throughout  the  Standard Conjectures, see [Grothendieck 1968]  \cite{Gro1}.) 
If   $V({\Bbb F}_q)$ is  defined  by a  system of polynomial equations,
 then the number of solutions of the system  is given by the formula:  
\begin{equation}\label{eq2}
|V({\Bbb F}_q)|=\sum_{i=0}^{2n}(-1)^i  ~tr~(Fr^i_q), 
\end{equation}
where  $tr$ is the trace of Frobenius  endomorphism  {\it loc. cit}.

 \bigskip 
  Let $B(V, {\cal L}, \sigma)$ be the twisted homogeneous coordinate ring of
an $n$-dimensional projective variety $V$ over a field $k$, where ${\cal L}$ is the invertible sheaf
of linear forms on $V$ and $\sigma$ an automorphism of $V$,  see
[Stafford \& van ~den ~Bergh 2001]  \cite{StaVdb1},
p. 180 for the notation and details.   Denote by ${\cal A}_V$ the {\it Serre $C^*$-algebra} 
of $V$,  i.e.  the norm-closure of a self-adjoint representation of   $B(V, {\cal L}, \sigma)$
by linear operators on a Hilbert space ${\cal H}$.  
  Consider a  stable $C^*$-algebra of ${\cal A}_V$,  i.e. the $C^*$-algebra  ${\cal A}_V\otimes {\cal K}$,
 where ${\cal K}$ is the $C^*$-algebra of compact operators on ${\cal H}$. 
 Let $\tau: {\cal A}_V\otimes {\cal K}\to {\Bbb R}$
be the unique normalized trace (tracial state) on  ${\cal A}_V\otimes {\cal K}$,  i.e. a positive linear functional
of norm $1$  such that $\tau(yx)=\tau(xy)$ for all $x,y\in {\cal A}_V\otimes {\cal K}$,  see  
 [Blackadar 1986] \cite{B},  p. 31.

Recall that ${\cal A}_V$ is the crossed product $C^*$-algebra of the form
${\cal A}_V\cong C(V)\rtimes {\Bbb Z}$,  where $C(V)$ is the 
commutative $C^*$-algebra of complex valued functions on $V$ 
and the product is taken by an automorphism of algebra $C(V)$ 
induced by the map $\sigma: V\to V$  [Nikolaev 2012] [arXiv:1208.2049].  
From the Pimsner-Voiculescu six term exact sequence for
crossed products,  one gets the  short exact sequence of algebraic $K$-groups:  
$0\to K_0(C(V))\buildrel  i_*\over\to  K_0({\cal A}_V)\to K_1(C(V))\to 0$, 
where   map  $i_*$  is induced by an  embedding of $C(V)$ 
into ${\cal A}_V$,   see   [Blackadar 1986]  \cite{B}, p. 83 for the details.  
 We  have $K_0(C(V))\cong K^0(V)$ and 
$K_1(C(V))\cong K^{-1}(V)$,  where $K^0$ and $K^{-1}$  are  the topological
$K$-groups of variety $V$, see  [Blackadar 1986]  \cite{B}, p. 80. 
By  the Chern character formula,  one gets
$K^0(V)\otimes {\Bbb Q}\cong H^{even}(V; {\Bbb Q})$ and 
$K^{-1}(V)\otimes {\Bbb Q}\cong H^{odd}(V; {\Bbb Q})$, 
where $H^{even}$  ($H^{odd}$)  is the direct sum of even (odd, resp.) 
cohomology groups of $V$.  
(Notice that $K_0({\cal A}_V\otimes {\cal K})\cong K_0({\cal A}_V)$ because
of  stability of the $K_0$-group with respect to tensor products by the algebra 
${\cal K}$,  see e.g.   [Blackadar 1986]  \cite{B}, p. 32.)
Thus one gets the following  commutative diagram:

\begin{picture}(300,100)(0,5)
\put(160,72){\vector(0,-1){35}}
\put(80,65){\vector(2,-1){45}}
\put(240,65){\vector(-2,-1){45}}
\put(10,80){$ H^{even}(V)\otimes {\Bbb Q} 
\buildrel  i_*\over\longrightarrow  K_0({\cal A}_V\otimes{\cal K})\otimes {\Bbb Q} 
\longrightarrow H^{odd}(V)\otimes {\Bbb Q}$}
\put(167,55){$\tau_*$}
\put(157,20){${\Bbb R}$}
\end{picture}

\noindent
where $\tau_*$ is a homomorphism  induced on $K_0$ by  the canonical  trace 
$\tau$ on the $C^*$-algebra  ${\cal A}_V\otimes {\cal K}$. 
Because   $H^{even}(V):=\oplus_{i=0}^n H^{2i}(V)$ and  
$H^{odd}(V):=\oplus_{i=1}^n H^{2i-1}(V)$,   one gets  for each  $0\le i\le 2n$ 
 an injective  homomorphism 
 \begin{equation}\label{eq2.6}
 H^i(V)\to  {\Bbb R}
 \end{equation}
 and we shall denote by $\Lambda_i$  an  additive abelian subgroup of real numbers
 defined by the homomorphism.    The $\Lambda_i$ is  known as   a {\it pseudo-lattice},
 see   [Manin 2004]  \cite{Man1},  Section 1.

Recall that  endomorphisms  of a pseudo-lattice are given as 
multiplication of points of $\Lambda_i$ by the real numbers $\alpha$
such that $\alpha\Lambda_i\subseteq\Lambda_i$.   It is known that
$End~(\Lambda_i)\cong {\Bbb Z}$ or $End~(\Lambda_i)\otimes {\Bbb Q}$
is a real algebraic number field such that $\Lambda_i\subset  End~(\Lambda_i)\otimes {\Bbb Q}$,
see e.g.  [Manin 2004]  \cite{Man1}, Lemma 1.1.1 for the case of quadratic fields. 
We shall write $\varepsilon_i$ to denote the unit of the order in the field $K_i:=End~(\Lambda_i)\otimes {\Bbb Q}$,  
which induces the  shift automorphism  of $\Lambda_i$,  
 see [Effros 1981]  \cite{E},  p. 38   for the details and terminology.

 Let $p$ be a good prime  and $V({\Bbb F}_q)$ a reduction
 of $V$ modulo $q=p^r$.   Consider a sub-lattice $\Lambda_i^{q}$ of $\Lambda_i$ of the index $q$; 
  by an  index of  the sub-lattice we understand  its  index as an abelian subgroup of $\Lambda_i$.
We shall write  $\pi_i(q)$ to  denote  an
 integer,   such that  multiplication by $\varepsilon_i^{\pi_i(q)}$ 
 induces  the shift automorphism of   $\Lambda_i^q$. 
  The trace of an algebraic number will be written as $tr~(\bullet)$.    
  Our main result relates invariants $\varepsilon_i$ and $\pi_i(q)$ of the $C^*$-algebra
 ${\cal A}_V$ to the cardinality of the set $V({\Bbb F}_q)$.   
\begin{thm}\label{thm1}
\qquad$|V({\Bbb F}_q)|=\sum_{i=0}^{2n}(-1)^i  ~tr~\left(\varepsilon_i^{\pi_i(q)}\right)$.  
\end{thm}
The article is organized as follows. Theorem \ref{thm1} is proved in Section 2.
Some examples can be found in Section 3.

\section{Proof of theorem \ref{thm1} }
We  shall split the proof in a series of lemmas starting with the following well-known
\begin{lem}\label{lm1}
There exists a symplectic unitary matrix $\Theta_q^i\in Sp ~(deg~P_i; ~{\Bbb R})$,  such that:
\begin{equation}\label{eq3}
Fr^i_q=q^{i\over 2}\Theta_q^i.  
\end{equation}
\end{lem}
{\it Proof.}  Recall that the eigenvalues of $Fr_q^i$ have absolute value $q^{i\over 2}$;
they come in the complex conjugate pairs.  On the other hand,  symplectic unitary matrices in
group   $Sp ~(deg~P_i; ~{\Bbb R})$ are known to have eigenvalues of absolute value $1$  
coming in complex conjugate pairs.   Since the spectrum of a matrix defines the  similarity class of 
matrix,  one can write the characteristic polynomial of $Fr_q^i$ in the form:
\begin{equation}\label{eq4}
 P_i(t)=det~(I-q^{i\over 2}\Theta_q^i t), 
\end{equation}
where matrix $\Theta_q^i\in Sp~(deg~P_i; ~{\Bbb Z})$ and its  eigenvalues have  absolute value $1$. 
It remains to compare (\ref{eq4}) with the formula:
\begin{equation}\label{eq5}
 P_i(t)=det~(I-Fr_q^i t), 
\end{equation}
i.e. $Fr_q^i=q^{i\over 2}\Theta_q^i$. Lemma \ref{lm1} follows.
$\square$

\begin{lem}\label{lm2}
Using a symplectic transformation one can bring matrix $\Theta_q^i$ to the block form:
\begin{equation}\label{eq6}
\Theta_q^i=\left(\matrix{A & I\cr -I & 0}\right), 
\end{equation}
where $A$ is a positive symmetric and $I$  the identity matrix.  
\end{lem}
{\it Proof.} Let  us write $\Theta_q^i$ in the block form:
\begin{equation}\label{eq7}
\Theta_q^i=\left(\matrix{A & B\cr C & D}\right), 
\end{equation}
where matrices $A,B,C, D$ are invertible and their transpose $A^T,B^T,C^T, D^T$
satisfy the symplectic equations: 
\begin{equation}\label{eq8}
\left\{
\begin{array}{cc}
A^TD-C^TB &= I,\\
A^TC-C^TA &= 0,\\
B^TD-D^TB &= 0. 
\end{array}
\right.
\end{equation}
Recall that symplectic matrices correspond to the linear
fractional  transformations  $\tau\mapsto {A\tau+B\over C\tau +D}$
of the Siegel half-space ${\Bbb H}_n=\{\tau=(\tau_j)\in {\Bbb C}^{{n(n+1)\over 2}}
~|~\Im (\tau_j)>0\}$ consisting of symmetric $n\times n$ matrices, see e.g. [Mumford 1983]
\cite{M},   p. 173.   One can always multiply    the nominator and denominator
of such a transformation  by $B^{-1}$  without affecting the transformation;
thus with no loss of generality,  we can assume that $B=I$.

We shall consider the symplectic matrix $T$ and its inverse $T^{-1}$
given by the formulas:
\begin{equation}\label{eq9}
T=\left(\matrix{I & 0\cr D & I}\right) \quad\hbox{and}
\quad T^{-1}=\left(\matrix{I & 0\cr -D & I}\right). 
\end{equation}

It is verified directly, that  
\begin{equation}\label{eq10}
T^{-1}\Theta_q^i T=
\left(\matrix{I & 0\cr -D & I}\right)
\left(\matrix{A & I\cr C & D}\right)
\left(\matrix{I & 0\cr D & I}\right)
=\left(\matrix{A+D & I\cr C-DA & 0}\right). 
\end{equation}

The system of equations (\ref{eq8})  with $B=I$ implies the 
following two equations:
\begin{equation}\label{eq11}
A^TD-C^T = I
 \quad\hbox{and}\quad
D=D^T. 
\end{equation}

Applying transposition to the both parts  of the first equation of  (\ref{eq11}),  
 one gets  $(A^TD-C^T)^T = I^T$ and,  therefore, $D^TA-C=I$.    
But the second of (\ref{eq11}) says that $D^T=D$;  thus one arrives 
at the equation $DA-C=I$.  The latter gives us $C-DA=-I$, 
which we substitute  in (\ref{eq10}) and get (in a new notation) the
conclusion of lemma \ref{lm2}. 

Finally, the middle of equations (\ref{eq8}) with $C=-I$ implies $A=A^T$,
i.e. $A$ is a symmetric matrix.   Since the eigenvalues of symmetric matrix
are always real and in view of  $tr~(A)>0$ (because $tr~(Fr_q^i)>0$),  one concludes that 
$A$ is similar to a positive matrix, see e.g.  [Handelman 1981]  \cite{Han1}, Theorem 1.   
Lemma \ref{lm2} follows. 
$\square$

\begin{lem}\label{lm3}
The symplectic unitary transformation $\Theta_q^i$ of $H^i(V; {\Bbb Z})$
descends to  an automorphism of $\Lambda_i$  given by the matrix: 
\begin{equation}\label{eq12}
M_q^i=\left(\matrix{A & I\cr I & 0}\right).  
\end{equation}
\end{lem}
{\it Proof.}
Since $\Lambda_i\subset K_i$ there exists a basis of $\Lambda_i$
consisting of algebraic numbers;  denote by $(\mu_1,\dots,\mu_k;   ~\nu_1,\dots,\nu_k)$ a basis
of  $\Lambda_i$ consisting of positive algebraic numbers  $\mu_i>0$ and 
$\nu_i>0$.  Using the injective homomorphism $\tau_*$ given by (\ref{eq2.6}) 
one can descend  $\Theta^i_q$  given by  (\ref{eq6})  to an  automorphism of $\Lambda_i$   so that 
\begin{equation}\label{eq13}
\left(\matrix{\mu'\cr \nu'}\right)
=\left(\matrix{A & I\cr -I & 0}\right)
\left(\matrix{\mu\cr \nu}\right)=
\left(\matrix{A\mu+\nu\cr -\mu}\right),
\end{equation}
where $\mu=(\mu_1,\dots,\mu_k)$ and $\nu=(\nu_1,\dots,\nu_k)$.
Because vectors $\mu$ and $\nu$ consist of positive entries and 
$A$ is a positive matrix,  it is immediate that  $\mu'=A\mu+\nu>0$ while
$\nu'=-\mu<0$.

Notice that all automorphisms in the (Markov) category of pseudo-lattices 
come from multiplication of the basis vector 
$(\mu_1,\dots,\mu_k;   ~\nu_1,\dots,\nu_k)$ of $\Lambda_i$ 
by an algebraic unit $\lambda>0$ of field $K_i$;   in particular, 
any such an automorphism must be given by a   
non-negative matrix,  whose  Perron-Frobenius eigenvalue coincides with $\lambda$.
Thus for any automorphism of $\Lambda_i$ it must hold 
$\mu'>0$ and $\nu'>0$.

In view of (\ref{eq13}),  we shall  consider an automorphism of $\Lambda_i$
given by  matrix $M_q^i=(A, I, I, 0)$;   clearly,  for $M_q^i$ it holds
 $\mu'=A\mu+\nu>0$ and  $\nu'=\mu>0$.   
Therefore  $M_q^i$ is a non-negative matrix satisfying  
  the   necessary condition  to belong to the  Markov category.  
   It is also a sufficient one,  because the similarity class of $M_q^i$ 
contains a representative whose Perron-Frobenius eigenvector can be taken for a
basis $(\mu, \nu)$ of $\Lambda_i$.  This argument finishes the proof of lemma \ref{lm3}.
$\square$

\begin{cor}\label{cr1}
$tr~(M_q^i)=tr~(\Theta_q^i)$. 
\end{cor}
{\it Proof.}  This fact is an implication of formulas (\ref{eq6}) and (\ref{eq12})
and  a direct computation  $tr~(M_q^i)=tr~(A)=tr~(\Theta_q^i)$.
$\square$

\begin{dfn}\label{df1}
We shall call $q^{i\over 2}M_q^i$ a Markov endomorphism of $\Lambda_i$
and denote it by $Mk_q^i$.  
\end{dfn}
\begin{lem}\label{lm4}
$tr~(Mk_q^i)=tr~(Fr_q^i)$.
\end{lem}
{\it Proof.}
Corollary \ref{cr1} says that $tr~(M_q^i)=tr~(\Theta_q^i)$,  and  
therefore:
\begin{equation}\label{eq14}
\begin{array}{ccc}
tr~(Mk_q^i) &=&  tr~(q^{i\over 2}M_q^i)= q^{i\over 2}~tr~(M_q^i)=\\
                     &=& q^{i\over 2} ~tr~(\Theta_q^i) =  tr~(q^{i\over 2}\Theta_q^i)=tr~(Fr_q^i). 
\end{array}
\end{equation}
In words,    Frobenius and Markov endomorphisms have the same trace,
i.e.  $tr~(Mk_q^i)=tr~(Fr_q^i)$.
$\square$

\medskip
\begin{rmk}\label{rm1}
{\normalfont
Notice  that,  unless $i$ or $r$  are even,  neither $\Theta_q^i$ nor $M_q^i$ are 
integer matrices;   yet    $Fr_q^i$ and $Mk_q^i$ are always  integer matrices.    
 }
\end{rmk}

\medskip
\begin{lem}\label{lm5}
There exists an algebraic unit $\omega_i\in K_i$ such that:

\medskip
(i) $\omega_i$ corresponds to the shift automorphism of an index
$q$ sub-lattice of pseudo-lattice $\Lambda_i$; 

\smallskip
(ii)  $tr~(\omega_i)=tr~(Mk_q^i)$.
\end{lem}
{\it Proof.}
To prove lemma \ref{lm5},  we shall  use the notion of a stationary dimension group and 
the corresponding shift  automorphism;  we refer the reader to  [Effros 1981]  \cite{E}, p. 37 and 
[Handelman  1981]  \cite{Han1}, p.57 for the notation and details  on  stationary dimension groups
and a survey of [Wagoner 1999]   \cite{Wag1}  for the general theory of subshifts of finite type.

Consider a stationary dimension group, $G(Mk_q^i)$, generated by  the Markov
endomorphism $Mk_q^i$:    
\begin{equation}\label{eq16}
{\Bbb Z}^{b_i} \buildrel Mk_q^i \over\to 
{\Bbb Z}^{b_i} \buildrel Mk_q^i \over\to 
{\Bbb Z}^{b_i} \buildrel Mk_q^i \over\to
\dots,
\end{equation}
where $b_i=deg~P_i(t)$.  
Let $\lambda_M$ be the Perron-Frobenius eigenvalue of matrix $M_q^i$. 
It is known, that $G(Mk_q^i)$ is  order-isomorphic to a  dense additive  abelian subgroup 
${\Bbb Z}[{1\over\lambda_M}]$ of  ${\Bbb R}$;  
here ${\Bbb Z}[x]$ is the set of all polynomials in one variable with the integer
coefficients.

Let $\widehat{Mk_q^i}$  be a  shift automorphism of  $G(Mk_q^i)$ [Effros 1981]  \cite{E},
p. 37.    To calculate the automorphism,  notice that multiplication 
of ${\Bbb Z}[{1\over\lambda_M}]$ by $\lambda_M$  induces an automorphism
of dimension group ${\Bbb Z}[{1\over\lambda_M}]$.   Since the determinant of matrix $M_q^i$ (i.e. the degree of 
Markov endomorphism)   is equal to $q^n$,   one concludes that  such an  automorphism corresponds to  a unit 
of the  endomorphism ring  of  a  sub-lattice of  $\Lambda_i$ of index $q^n$.
We shall denote such a unit by $\omega_i$.    Clearly,  $\omega_i$ generates 
the required shift automorphism $\widehat{Mk_q^i}$  through  multiplication 
of  dimension group  ${\Bbb Z}[{1\over\lambda_M}]$  by the algebraic number  $\omega_i$. 
Item (i) of lemma \ref{lm5} follows. 

\medskip
Consider the Artin-Mazur zeta function of $Mk_q^i$:
\begin{equation}\label{eq17}
\zeta_{Mk_q^i}(t)=\exp\left(\sum_{k=1}^{\infty}{tr~\left[(Mk_q^i)^k\right]\over k}t^k\right)
\end{equation}
and such of $\widehat{Mk_q^i}$:
\begin{equation}\label{eq18}
\zeta_{\widehat{Mk_q^i}}(t)=\exp\left(\sum_{k=1}^{\infty}{tr~\left[(\widehat{Mk_q^i})^k\right]\over k}t^k\right).
\end{equation}
Since  $Mk_q^i$ and $\widehat{Mk_q^i}$ are shift equivalent matrices,  one concludes 
 that  $\zeta_{Mk_q^i}(t)\equiv \zeta_{\widehat{Mk_q^i}}(t)$,    see [Wagoner 1999]  \cite{Wag1}, p. 273.    
 In particular, 
\begin{equation}\label{eq19} 
tr~(Mk_q^i)=tr~(\widehat{Mk_q^i}).
\end{equation}
But   $tr~(\widehat{Mk_q^i})=tr~(\omega_i)$,  where on the right hand side is the
trace of an algebraic number.  In view of (\ref{eq19}) one gets the conclusion of 
item (ii)  of lemma \ref{lm5}. 
$\square$

\begin{lem}\label{lm6}
There exists a positive integer $\pi_i(q)$,  such that:
\begin{equation}\label{eq19.5} 
\omega_i=\varepsilon_i^{\pi_i(q)},
\end{equation}
where $\varepsilon_i\in End~(\Lambda_i)$ is the fundamental unit  
corresponding to  the shift automorphism of pseudo-lattice $\Lambda_i$. 
\end{lem}
{\it Proof.}  
Given an automorphism $\omega_i$ of a finite-index sub-lattice of $\Lambda_i$
one can extend $\omega_i$ to an automorphism of entire $\Lambda_i$,
since $\omega_i\Lambda_i=\Lambda_i$.  Therefore each unit of (endomorphism ring of)
a sub-lattice is also a unit of the host pseudo-lattice.  Notice that the converse statement 
is false in general. 

By virtue of the Dirichlet Unit Theorem  each  unit of $End~(\Lambda_i)$ 
is a product of a finite number of (powers of)  fundamental units of   $End~(\Lambda_i)$.
We shall denote by $\pi_i(q)$ the least positive integer,  such that  $\varepsilon_i^{\pi_i(q)}$
is the shift automorphism of a sub-lattice of index $q$ of pseudo-lattice $\Lambda_i$. 
The number $\pi_i(q)$ exists and uniquely defined,  albeit no general formula for its calculation 
 is known,   see remark \ref{rm2}.   It is clear from construction,  that $\pi_i(q)$ 
 satisfies the claim of lemma \ref{lm6}.
 $\square$

\medskip
\begin{rmk}\label{rm2}
{\normalfont
No general formula for the number $\pi_i(q)$ as a function of $q$ 
is known;  however,  if the rank of $\Lambda_i$ is two (i.e. $n=1$), 
then there are classical results recorded in  e.g.   [Hasse 1950]  \cite{HS},  p.298. 
See also examples section of this paper.   
}
\end{rmk}

\bigskip
Theorem \ref{thm1} follows from formula (\ref{eq2}) and lemmas \ref{lm4}-\ref{lm6}.
$\square$

\section{Examples}
We shall consider two sets  of examples illustrating theorem \ref{thm1};  both
deal with  non-singular elliptic curves defined over the field of  complex numbers. 
We shall assume that ${\cal E}_{\tau}$ is such a curve isomorphic to a complex torus,
i.e.  ${\cal E}_{\tau}={\Bbb C}/({\Bbb Z}+{\Bbb Z}\tau)$,
where $\tau\in {\Bbb H}=\{x+iy\in {\Bbb C} ~|~y>0\}$ is  a complex modulus, 
see e.g. [Hartshorne  1977]  \cite{H},  p. 326.
The Serre $C^*$-algebra ${\cal A}_{\cal E_{\tau}}$  of elliptic curve ${\cal E}_{\tau}$ is known
to be isomorphic to the so-called  noncommutative torus ${\cal A}_{\theta}$ 
with  the  unit scaled by  a constant $0<\log\mu<\infty$,   where  $\theta$ is irrational  
and $\mu$ a  positive real number  [Nikolaev 2011]  [arXiv:1109.6688];
we refer the reader to [Rieffel 1990]  \cite{Rie1}  for the definition and properties 
of noncommutative tori. 
It is known that $K_0({\cal A}_{\theta})=K_1({\cal A}_{\theta})\cong {\Bbb Z}^2$
so that the canonical trace $\tau$ on ${\cal A}_{\theta}$ provides us with the
following useful formula: 
\begin{equation}\label{eq20} 
\tau_*(K_0({\cal A}_{\cal E_{\tau}}\otimes {\cal K}))=\mu({\Bbb Z}+{\Bbb Z}\theta).
\end{equation}
 Because  $H^0({\cal E}_{\tau}; {\Bbb Z})=H^2({\cal E}_{\tau};  {\Bbb Z})\cong {\Bbb Z}$ while  
 $H^1({\cal E}_{\tau};  {\Bbb Z})\cong {\Bbb Z}^2$,  one gets 
 from  injective homomorphism (\ref{eq2.6})   the following collection of
 pseudo-lattices: 
\begin{equation}\label{eq21} 
\Lambda_0=\Lambda_2\cong {\Bbb Z} \quad\hbox{and}\quad  
\Lambda_1\cong\mu({\Bbb Z}+{\Bbb Z}\theta).
\end{equation}

\subsection{Complex multiplication}
To get an arithmetic variety,  we shall assume  that ${\cal E}_{\tau}$ has 
complex multiplication;  the multiplication is characterized by the endomorphism
ring of ${\cal E}_{\tau}$   which  is  an order of conductor
$f\ge 1$ in the imaginary quadratic field ${\Bbb Q}(\sqrt{-D})$,  
see e.g. [Hartshorne  1977]  \cite{H},  p. 330.   In this case $\tau\in {\Bbb Q}(\sqrt{-D})$
  ${\cal E}_{\tau}\cong {\cal E}(k)$,  where $k$ is the Hilbert
class field (i.e. maximal abelian extension) of the field ${\Bbb Q}(\sqrt{-D})$.    
 For such a   curve   formulas (\ref{eq21})   will depend on $f$ and $D$:  
\begin{equation}\label{eq22} 
\Lambda_0=\Lambda_2\cong {\Bbb Z} \quad\hbox{and}\quad
\Lambda_1 = \varepsilon[{\Bbb Z}+(f\omega){\Bbb Z}], 
\end{equation}
where  $\omega={1\over 2}(1+\sqrt{D})$ if $D\equiv 1~mod~4$ and $D\ne 1$  or $\omega=\sqrt{D}$ if $D\equiv 2,3~mod~4$
and  $\varepsilon>1$ is the fundamental unit of  order ${\Bbb Z}+(f\omega){\Bbb Z}$.    
As expected, $\Lambda_1\subset K_1$,   where $K_1$  is  the real  quadratic  field
${\Bbb Q}(\sqrt{D})$.

Let $p$ be a good prime.  Consider a localization ${\cal E}({\Bbb F}_p)$
of curve  ${\cal E}(k)$ at the prime ideal ${\goth P}$ over $p$. 
It is well known, that the Frobenius endomorphism  of elliptic curve
with complex multiplication  is defined  by the so-called  Gr\"ossencharacter,
which is essentially  a complex number  $\alpha_{{\goth P}}\in {\Bbb Q}(\sqrt{-D})$
of absolute value $\sqrt{p}$;
multiplication of the lattice $L_{CM}={\Bbb Z}+{\Bbb Z}\tau$ by $\alpha_{{\goth P}}$
induces the Frobenius endomorphism $Fr_p^1$ on $H^1({\cal E}(k); {\Bbb Z})$,
see e.g.  [Silverman 1994] \cite{S},  p. 174.  
Thus   one arrives  at the following matrix form for the Frobenius \&  Markov endomorphisms 
and the shift automorphism,  respectively: 
\begin{equation}\label{eq23}
\left\{
\begin{array}{cc}
Fr_p^1 &=  \left(\matrix{tr~(\alpha_{{\goth P}}) & p\cr -1 & 0}\right),\\
Mk_p^1 &=  \left(\matrix{tr~(\alpha_{{\goth P}}) & p\cr 1 & 0}\right), \\
\widehat{Mk_p^1} &=  \left(\matrix{tr~(\alpha_{{\goth P}}) & 1\cr  1 & 0}\right) . 
\end{array}
\right.
\end{equation}
To calculate  positive integer $\pi_1(p)$ appearing in theorem \ref{thm1},
  denote by $\left({D\over p}\right)$  the Legendre symbol of $D$ and $p$. 
    A classical result of the theory of real quadratic fields asserts that     
$\pi_1(p)$  must be one of the divisors of the integer number:
\begin{equation}\label{eq24} 
p-\left({D\over p}\right),
\end{equation}
see e.g. [Hasse  1950]  \cite{HS}, p. 298.   Thus the trace of Frobenius 
endomorphism on $H^1({\cal E}(p); {\Bbb Z})$ is given by the formula: 
\begin{equation}\label{eq24.5} 
tr~(\alpha_{{\goth P}})=tr~(\varepsilon^{\pi_1(p)}),
\end{equation}
where $\varepsilon$ is taken from (\ref{eq22}).  The right hand side of 
(\ref{eq24.5}) can be further simplified,   since 
\begin{equation}\label{eq24.6} 
tr~(\varepsilon^{\pi_1(p)})= 2T_{\pi_1(p)}\left[~{1\over 2}~ tr~(\varepsilon)\right],
\end{equation}
where $T_{\pi_1(p)}(x)$ is the Chebyshev polynomial (of the first kind) 
of degree $\pi_1(p)$.  Thus  one obtains a  formula for the number of (projective) solutions
of a cubic equation over field  ${\Bbb F}_p$  in terms of invariants of
pseudo-lattice $\Lambda_1$:         
\begin{equation}\label{eq25} 
|{\cal E}({\Bbb F}_p)|=1 + p - 2T_{\pi_1(p)}\left[ ~{1\over 2}~ tr~(\varepsilon)\right] . 
\end{equation}

\subsection{Rational elliptic curve}
Suppose that  $b\ge 3$ is an integer and  consider a rational elliptic curve in 
complex projective  plane given by a homogeneous Legendre equation:
\begin{equation}\label{eq26} 
y^2z=x(x-z)\left(x-{b-2\over b+2}z\right).  
\end{equation}
The Serre $C^*$-algebra of projective variety (\ref{eq26}) is isomorphic
(modulo an ideal) to the so-called Cuntz-Krieger algebra ${\cal O}_B$, where 
\begin{equation}\label{eq27} 
B=\left(\matrix{b-1 & 1\cr b-2 & 1}\right)  
\end{equation}
is a positive integer matrix [Nikolaev 2012]  [arXiv:1201.1047];  
for the definition and properties of algebra ${\cal O}_B$ 
we refer the reader to [Cuntz \& Krieger 1980] \cite{CuKr1}.

Recall that ${\cal O}_B\otimes {\cal K}$ is the   crossed product $C^*$-algebra
of a stationary AF $C^*$-algebra by its shift automorphism, 
see [Blackadar 1986]  \cite{B},  p. 104; 
the AF $C^*$-algebra has the following dimension group:
\begin{equation}\label{eq28}
{\Bbb Z}^2 \buildrel B^T \over\to 
{\Bbb Z}^2\buildrel B^T \over\to 
{\Bbb Z}^2 \buildrel B^T\over\to
\dots,
\end{equation}
where $B^T$ is the transpose of matrix $B$.   
Because $\mu$ in formula (\ref{eq21}) must be
a positive eigenvalue of matrix $B^T$,   one gets:   
\begin{equation}\label{eq30a}
\mu  =  {2-b+\sqrt{b^2-4}\over 2}. 
\end{equation}
Likewise,  since $\theta$ in formula (\ref{eq21}) must be the corresponding 
positive eigenvector $(1,\theta)$ of the same matrix, one gets:
\begin{equation}\label{eq30b}
\theta =   {1\over 2}\left(\sqrt{{b+2\over b-2}}-1\right).
\end{equation}
Therefore,  pseudo-lattices $\Lambda_i$ are  $\Lambda_0=\Lambda_2\cong {\Bbb Z}$
and 
\begin{equation}\label{eq29} 
\Lambda_1\cong 
{2-b+\sqrt{b^2-4}\over 2} 
\left[{\Bbb Z}+
{1\over 2}\left(\sqrt{{b+2\over b-2}}-1\right)
{\Bbb Z}\right].
\end{equation}
As expected,  pseudo-lattice $\Lambda_1\subset K_1$,  where $K_1={\Bbb Q}(\sqrt{b^2-4})$
is a real quadratic field.  

\bigskip
Let $p$ be a good prime and let ${\cal E}({\Bbb F}_p)$ be the reduction
of  curve (\ref{eq26})  modulo  $p$.  Using the same argument as in Section 3.1,
we determine number $\pi_1(p)$ as one of the divisors of integer number: 
\begin{equation}\label{eq32} 
p-\left({b^2-4\over p}\right).
\end{equation}
Unlike the case of complex multiplication, 
the Gr\"ossencharacter is no longer  available for elliptic curve   (\ref{eq26});   
yet  the trace of Frobenius endomorphism  can be computed using 
theorem \ref{thm1}: 
\begin{equation}\label{eq31} 
tr~(Fr_p^1)=tr~\left[(B^T)^{\pi_1(p)}\right].
\end{equation}
Using the Chebyshev polynomials,  one can write (\ref{eq31}) in the 
form:
\begin{equation}\label{eq31.5} 
tr~(Fr_p^1)=2T_{\pi_1(p)}\left[~{1\over 2}~ tr~(B^T)\right]. 
\end{equation}
In view of (\ref{eq27}) one gets $tr~(B^T)=b$,  so that (\ref{eq31.5}) takes the form:
\begin{equation}\label{eq31.6} 
tr~(Fr_p^1)=2T_{\pi_1(p)}\left({b\over 2}\right). 
\end{equation}
Thus  one obtains a  formula for the number of  solutions
of equation (\ref{eq26}) over field  ${\Bbb F}_p$  in terms of invariants of
pseudo-lattice $\Lambda_1$:         
\begin{equation}\label{eq33} 
|{\cal E}({\Bbb F}_p)|=1 + p - 2T_{\pi_1(p)}\left({b\over 2}\right).
\end{equation}
We conclude by an example comparing formula (\ref{eq33}) with
the known results for a rational elliptic curve  in the Legendre form, 
see e.g. [Hartshorne  1977]  \cite{H},  p.  333  and  [Kirwan 1992]  \cite{K}, 
pp. 49-50.

\bigskip
\begin{exm}\label{ex1}
{\normalfont
Suppose that $b\equiv 2~mod~4$.   Recall that the $j$-invariant takes the same value on 
$\lambda$, $1-\lambda$ and ${1\over\lambda}$,   see   e.g.  [Hartshorne  1977]  \cite{H},  p.  320.
Therefore,   one can bring (\ref{eq26}) to the form:
\begin{equation}\label{eq34} 
y^2z=x(x-z)(x-\lambda z),    
\end{equation}
 where $\lambda={1\over 4}(b+2)= 2,3,4,\dots$     Notice that for curve (\ref{eq34}):
\begin{equation}\label{eq35} 
tr~(B^T)=b=2(2\lambda-1).     
\end{equation}
 To calculate (\ref{eq33}) for elliptic curve (\ref{eq34}),   recall that in view of (\ref{eq31.6})
 one gets:
\begin{equation}\label{eq37} 
  tr~(Fr_p^1)=2 ~T_{\pi_1(p)} (2\lambda-1).  
\end{equation}
It will be useful to express  Chebyshev polynomial   $T_{\pi_1(p)} (2\lambda-1)$   
in terms of the hypergeometric function $_2F_1(a, b; c; z)$;  the standard 
formula brings (\ref{eq37}) to the form: 
\begin{equation}\label{eq38} 
  tr~(Fr_p^1)=2 ~_2F_1(-\pi_1(p), ~\pi_1(p); ~{1\over 2}; ~1-\lambda).  
\end{equation}
We leave  to the reader to prove the identity:
\begin{equation}\label{eq39}
\begin{array}{cc}
  2 ~_2F_1(-\pi_1(p), ~\pi_1(p); ~{1\over 2}; ~1-\lambda) =&\\
   &\\
=  (-1)^{\pi_1(p)} ~_2F_1(\pi_1(p)+1,~\pi_1(p)+1; ~1; ~\lambda). &
\end{array}
\end{equation}
In the last formula: 
\begin{equation}\label{eq40} 
 _2F_1(\pi_1(p)+1,~\pi_1(p)+1; ~1; ~\lambda)=
 \sum_{r=0}^{\pi_1(p)}\left(\matrix{\pi_1(p)\cr r}\right)^2 \lambda^r,  
\end{equation}
see  [Carlitz 1966]  \cite{Car1},  p.328.

 Recall that $\pi_1(p)$ is a divisor of (\ref{eq32}), which in our case 
 takes the value ${p-1\over 2}$.  Bringing together formulas (\ref{eq33}),
 (\ref{eq38})-(\ref{eq40}) one gets:
\begin{equation}\label{eq41} 
|{\cal E}({\Bbb F}_p)|=1+p + (-1)^{{p-1\over 2}} 
\sum_{r=0}^{{p-1\over 2}} \left(\matrix{{p-1\over 2} \cr r}\right)^2 \lambda^r, 
\end{equation}
compare with  [Kirwan 1992]  \cite{K},   pp. 49-50 and 
 [Hartshorne  1977]  \cite{H},  p.  333 for a relation with the
 Hasse  invariant.  
 
  }
\end{exm}



\vskip1cm

\textsc{The Fields Institute for Research in Mathematical Sciences, Toronto, ON, Canada,  
E-mail:} {\sf igor.v.nikolaev@gmail.com}

\smallskip
{\it Current address: 1505-657 Worcester St.,  Southbridge,  MA 01550,  U.S.A.}


\begin{thebibliography}{100}
\bibitem{B}
B.~Blackadar, $K$-Theory for Operator Algebras, MSRI Publications,
Springer, 1986

\bibitem{Car1}
L.~Carlitz,  Some binomial coefficient identities,  Fibonacci Quart. 4 (1966), 
323-331.  


\bibitem{CuKr1}
J. ~Cuntz and W. ~Krieger, A class of $C^*$-algebras and topological Markov
chains, Invent. Math. 56 (1980), 251-268. 


\bibitem{E}    
E.~G.~Effros, Dimensions and $C^*$-Algebras, in: Conf. Board of the Math.
Sciences, Regional conference series in Math., No.46, AMS,  1981.



\bibitem{Gro1}
A.~Grothendieck,  Standard conjectures on algebraic cycles,
in:  Algebraic Geometry,  Internat. Colloq. Tata Inst. Fund. Res., Bombay, 1968.


\bibitem{Han1}
D.~Handelman, Positive matrices and dimension groups affiliated
to $C^*$-algebras and topological Markov chains, J. Operator
Theory 6 (1981), 55-74.


\bibitem{H}
R.~Hartshorne, Algebraic Geometry, GTM 52, Springer, 1977.  


\bibitem{HS}
H.~Hasse, Vorlesungen \"uber Zahlentheorie, Springer, 1950. 


\bibitem{K}
F.~Kirwan,  Complex Algebraic Curves,  LMS Student Texts 23,
Cambridge, 1992.  


\bibitem{Man1}
Yu.~I.~Manin, Real multiplication and noncommutative geometry,
in ``Legacy of Niels Hendrik Abel'', 685-727, Springer, 2004.


\bibitem{M}
D.~Mumford, Tata Lectures on Theta I, Birkh\"auser, 1983.   




\bibitem{Rie1}
M.~A.~Rieffel, Non-commutative tori -- a case study of non-commutative
differentiable manifolds,  Contemp. Math. 105 (1990), 191-211. 
Available {\sf http://math.berkeley.edu/$\sim$rieffel/}  



\bibitem{S}
J.~H.~Silverman, Advanced Topics in the Arithmetic of Elliptic Curves,
GTM 151, Springer 1994.



\bibitem{StaVdb1}
J.~T.~Stafford and M.~van ~den ~Bergh, Noncommutative curves and noncommutative
surfaces, Bull. Amer. Math. Soc. 38 (2001), 171-216. 





\bibitem{Wag1}
J.~B.~Wagoner, Strong shift equivalence theory and the shift equivalence
problem, Bull. Amer. Math. Soc. 36 (1999), 271-296. 





\end{thebibliography}
\end{document}